\documentstyle[12pt]{article}
\title{The First Eigenvalue of $P$-Manifolds}
\author{Akhil Ranjan and G.~Santhanam }
\date{September 7, 1995}
\newtheorem{theorem}{Theorem}
\newtheorem{lemma}{Lemma}
\newtheorem{propn}{Proposition}

\newtheorem{sublemma}{Sublemma}
\newcommand{\ncmd}{\newcommand}
\ncmd{\cm}{\mbox{${D_{\max}}$}}
\ncmd{\lpx}{\mbox{${\mu_{p(x)}}$}}
\ncmd{\lpxx}{\mbox{${\mu_{p(x)-1}}$}}
\ncmd{\wm}{\mbox{$\widetilde{M}$}}
\ncmd{\wg}{\mbox{$\widetilde{g}$}}
\ncmd{\wgaa}{\mbox{$\widetilde{\gamma_1}$}}
\ncmd{\wga}{\mbox{$\widetilde{\gamma}$}}
\ncmd{\rpn}{\mbox{$I\!\!RI\!\!P^n$}}
\ncmd{\minnu}{\mbox{$\parallel\!\nabla f\!\parallel$}}
\ncmd{\dt}{\mbox{${d\over{dt}}$}}
\ncmd{\dte}{\mbox{${{\partial}\over{\partial\theta}}$}}
\ncmd{\jt}{\mbox{${J_{\theta}}$}}
\ncmd{\euu}{\mbox{$E_{-{{1+\cos t}\over 2}}$}}
\ncmd{\eu}{\mbox{$E_{-\cos t}$}}
\ncmd{\eeu}{\mbox{$E_{{{1-\cos t}\over 2}}$}}
\ncmd{\ptp}{\mbox{$P_{2\pi}$}}
\ncmd{\orr}{\mbox{$\overline{R}$}}
\ncmd{\ste}{\mbox{$\sigma(\theta)$}}
\ncmd{\nabu}{\mbox{$\nabla f$}}
\ncmd{\nabx}{\mbox{$\nabla _X$}}
\ncmd{\navu}{\mbox{$\nabla_{{\nabla u}}$}}
\ncmd{\oga}{\mbox{$\overline{\gamma}$}}
\ncmd{\ogag}{\mbox{$\overline{\gamma_1}$}}
\ncmd{\olim}{\mbox{$\overline{m}$}}
\ncmd{\oliu}{\mbox{$\overline{u}$}}
\ncmd{\olimm}{\mbox{$\overline{M}$}}
\ncmd{\orc}{\mbox{$\overline{R}$}}
\ncmd{\ioga}{\mbox{$I_{\overline{\gamma}}$}}
\ncmd{\oho}{\mbox{$\overline{h}$}}
\ncmd{\grhu}{\mbox{$\nabla_hu$}}
\ncmd{\grhau}{\mbox{$\nabla_{\oho}u$}}
\begin{document}
\maketitle
\begin{abstract}
Antonio Ros gave a lower bound for the first eigenvalue $\lambda_1$ of
$\Delta$ of a $P$-manifold $(M, g)$ in terms of the lower bound on the
Ricci curvature $Ric_M$ and asked what happened when this lower bound was
achieved. In this paper we look in to this question and show that there
are strong implications on the geometry and topology of the underlying
manifold. In particular we show that in
case of spheres or real projective spaces we have isometry with the
standard metric. In other cases, with some additional hypothesis, we again
show isometry with standard models.
\end{abstract}
\footnote{1991 Mathematics Classification: 53C20, 53C22, 53C35,
58G25, 58G26.}
\section{Introduction}
Let $(M, g)$ be a compact Riemannian manifold, $\Delta$ the Laplacian of
$(M, g)$ and $Spec(M, g):=\{\lambda_1 < \lambda_2 <, \cdots \}$ the
spectrum of $\Delta$ of $(M, g)$.

It is an important problem in geometry to find lower bounds for the
eigenvalues of $\Delta$ of $(M, g)$ in terms of the given geometric data
and characterize those Riemannian manifolds $(M, g)$ for which these
lower bounds are attained. Lichnerowicz proved in \cite{L} that
{\it if $(M, g)$ is a complete Riemannian manifold of dimension $n\geq 2$
with Ricci curvature $Ric_M \geq l$, then
the first eigenvalue $\lambda_1$ satisfies the inequality
$\lambda_1\geq {n\over {n-1}}l$}.
Later Obata proved in \cite{MO} that {\it the equality is
attained only for the round sphere of radius ${1\over l}$.}
considering this problem for $P$-manifolds,
Antonio Ros proved in \cite{AR} that {\it if $(M, g)$ is a
$P_{2\pi}$-manifold of dimension $n \geq 2$ with Ricci curvature
$Ric_M\geq l$, then the
first eigenvalue $\lambda_1$ satisfies the
inequality $\lambda_1 \geq {1\over 3}(2l+n+2)$ and the equality is
attained iff for any first eigenfunction $f$ we have that
$f(\gamma_u(t))= A_u\cos t +B_u\sin t + C_u$ for $u\in UM$}. He further
remarked that in view of Obata's theorem this can happen only for a
small class of manifolds.

In this paper we substantiate his claim by proving
\begin{theorem}
Let $(M, g)$ be a $P_{2\pi}$-manifold of dimension $n\geq 2$ with
Ricci curvature $Ric_M\geq l$ and $\lambda_1={1\over 3}(2l+n+2)$.
Then
\begin{enumerate}
\item
\begin{enumerate}
\item $\lambda_1= {{k(m+1)} \over 2}=\lambda_1(\olimm)$
and $l=Ric_{\olimm}$ where $\olimm$
is a compact rank-1 symmetric space(CROSS) of dimension $n=km$ with
sectional curvature ${1\over 4}\leq K_{\olimm}\leq 1$ and $k=1$, $2$,
$4$, $8$ or $n$ is degree of the generator of
$H^*(M, I\!\!\!Q)=H^*(\olimm, I\!\!\!Q)$ and
$H^*(M, Z\!\!Z_2)=H^*(\olimm, Z\!\!Z_2)$.
\item If $k\geq 2$ then $M$ is simply connected and the integral
cohomology ring of $M$ is same as that of $\olimm$.
\end{enumerate}
\item If $k=1$ then $(M, g)$ is isometric to $\rpn$ with
constant sectional curvature ${1\over 4}$.
\item If $k=n$ then $(M, g)$ is isometric to $S^n$ with constant sectional
curvature $1$. ({\bf Lichnerowicz-Obata theorem})
\item If $k=2$, $4$ or $8$ and if there is a first eigenfunction $f$
without saddle points then $(M, g)$ is isometric to $\olimm$ of dimension
$km$.
\end{enumerate}
\end{theorem}
The main step in the proof of theorem 1 is the following
\begin{theorem}
Let $(M, g)$ be a $\ptp$-manifold of dimension $n\geq 2$ and $\lambda$ be
an eigenvalue of $\Delta$ with an eigenfunction $f$ such that
$f(\gamma_u(t)) = A_u\cos t + B_u\sin t + C_u$ for $u\in UM$. Then
$\lambda={{k(m+1)}\over 2}=\lambda_1(\olimm)$ where $\olimm$ is as in
theorem 1.
\end{theorem}
We refer to \cite{BE} and \cite{CE} for definitions, basic tools and
results used in this paper.
\section{Preliminaries}
In this section we study the topology of critical sets of the function
$f$ of the form $f(\gamma_u(t))= A_u\cos t + B_u\sin t + C_u$ for $u\in
UM$ on a $\ptp$-manifold $(M, g)$.

\noindent{\bf Definition:} Let $(M, g)$ a complete Riemannian manifold. A
subset $B \subseteq M$ is called {\it totally $a$-convex} if for any
pair of points $a_1$, $a_2\in B$ and any geodesic
$\gamma: [0, r]\to M$ with $\gamma(0)= a_1$ and $\gamma(r)=a_2$ and
$r < a$ then $\gamma([0, r])\subseteq B$. (See \cite{GG}).
\begin{theorem}
Let $(M, g)$ be a $\ptp$ manifold and $f\in C^{\infty}(M)$ be such that
$f(\gamma_u(t))= A_u\cos t + B_u\sin t + C_u$ for $u\in UM$. Then
\begin{enumerate}
\item For each critical value $\alpha$ of the function $f$, the
set $D_{\alpha}:=\{x\in M:f(x)=\alpha ~\&~ \nabla f(x)=0\}$ is a
totally $2\pi$-convex, totally geodesic submanifold of $(M, g)$.
\item $d(D_{\alpha}, D_{\beta})=\pi$ for $\alpha\neq \beta$.
\item The function $f$ has only finitely many critical values.
\end{enumerate}
\end{theorem}
\noindent{\bf Proof of theorem 3:}
Let $x\in M$. Then
$f(\gamma_u(t))=A_u\cos t + B_u\sin t + C_u$ for every  $u\in U_xM$,
the unit sphere in $T_xM$. If $x$ is a critical point of the function
$f$, then, since $\nabla f(x)=0$, we have that
\begin{eqnarray*}
B_u & = & \dt\mid_{t=0}f(\gamma_u(0)) \\
{}  & = & <\nabla f(x), \gamma_u'(0)> \\
{}  & = & 0
\end{eqnarray*}
Therefore $f(\gamma_u(t))=A_u\cos t + C_u$ for every $u\in U_xM$ and
$x$ a critical point of the function $f$.

\noindent{\bf Proof of 1:} Let $x$, $y\in D_{\alpha}$ and $\gamma_u$ be a
geodesic joining $x$ and $y$ such that $\gamma_u(0)=x$ and
$\gamma_u(r)=y$ for some $r\in I\!\!R^+$. Since $f(x)=f(y)=\alpha$ and
$f(\gamma_u(t))=A_u\cos t + C_u$, we have that
$A_u + C_u=A_u\cos r + C_u$. Hence $A_u=0$ if $r<2\pi$.
This shows that $f(\gamma_u(t))=\alpha$ for all $t\in [0, r]$ and hence
$\gamma([0, r])\subseteq D_{\alpha}$.
This proves that $D_{\alpha}$ is totally $2\pi$ convex.

To show that $D_{\alpha}$ is totally geodesic let us start with a
$u\in UD_{\alpha}$, the unit tangent bundle of $D_{\alpha}$.
Then, since $0=\nabla^2f(u, u)=-A_u$, we have that
$f(\gamma_u(t))=\alpha$,
a constant. Hence $\gamma_u\subseteq D_{\alpha}$. Since $(M, g)$ is a
$\ptp$-manifold, this remains true even if $u$ is a unit tangent vector
based at a boundary point of $D_{\alpha}$. Hence $D_{\alpha}$ is a
totally geodesic submanifold of $(M, g)$.

\noindent{\bf Proof of 2:} Let $\alpha$ and $\beta$ be two
critical values of the function $f$ such that $\alpha\neq \beta$.
Let $x\in D_{\alpha}$, $y\in D_{\beta}$ with $d(x, y)=t_0$
for some $t_0\in I\!\!R^+$ and $\gamma_u$ be a geodesic
segment such that $\gamma_u(0)=x$ and $\gamma_u(t_0)=y$.
Then $f(\gamma_u(t))= A_u\cos t + C_u$ and
\begin{eqnarray*}
-A_u\sin t_0 & = & \dt\mid_{t=t_{0}} f(\gamma_u(t)) \\
{}           & = & < \nabla f(y), \gamma_u'(t_0) > \\
{}           & = & 0
\end{eqnarray*}
This can happen only if $t_0=\pi$.
This proves that $d(D_{\alpha}, D_{\beta}) = \pi$ for $\alpha\neq\beta$.

\noindent{\bf Proof of 3:} It is obvious as the critical submanifolds are
constant distance apart.
\subsection{}
Since the function $f$ has only finitely many critical values, we denote
these critical values by $\max(f)=\alpha_1, \alpha_2, \cdots ,
\alpha_p=\min(f)$ and we denote by $D_i$ the critical submanifold
$\{x\in M: f(x)=\alpha_i  \&  \nabla f(x)=0\}$.

Let $x\in D_{\max}=\{ x\in M:f(x)=\max(f)\}$.
Then $- \nabla^2f(x)$ is positive semi-definite for each $x\in \cm$.
Therefore we can write the eigenvalues
of $- \nabla^2 f(x)$ as $\lpx > \lpxx > \cdots >\mu_2(x) >\mu_1=0$
where each $\mu_i(x)$ is a function on $\cm$ for $1\leq i \leq p(x)$ and
$p(x)\in \{1, 2, \cdots , n\}$.

For each $i$, we denote by $E_{\mu_{i(x)}}$, the $\mu_i(x)$-
eigensubspace of $-\nabla^2f(x)$, by $U_{\mu_{i(x)}}$ the unit sphere in
$E_{\mu_{i(x)}}$ and by $S_{\mu_{i(x)}}(0, r)$ the sphere of radius
$r$ in $E_{\mu_{i(x)}}$.

Let $u\in U{_{\mu_i(x)}}$. Then $\max(f)=A_u + C_u$ and
$\mu_i(x) = -\nabla^2f(u, u) = A_u$. Therefore $A_u$ and hence
$C_u = \max(f)- A_u$ are constants on $U{_{\mu_i(x)}}$. Now we define
$S(\mu_i(x), r) :=\exp_x(S_{\mu_{i(x)}}(0, r))$, the exponential
image of the sphere $S_{\mu_{i(x)}}(0, r)$ of radius $r$. Since the
function $f$ is constant on $S(\mu_i(x), r)$ for each $r$, we have that
$\nabla f$ is normal to $S(\mu_i(x), r)$. Therefore the geodesics
$\gamma_u$ are integral curves of $\nabla f$ for $u\in U_{\mu_{i(x)}}$
and $\nabla f = - \mu_i(x) \sin t \partial_t$. From this it follows that
$\nabla f(y)=0$ for $y\in D_i(x) := S(\mu_i(x), \pi)$.
Further $D_i(x) =\{ y\in M:f(y)=\max(f)-2\mu_i(x) ~\&~ \nabla f(y)=0\}$.
This can be seen as follows: Let $y\in D_i(x)$.
Then
\begin{eqnarray*}
f(y) & = & f(\gamma_u(\pi)) \\
{}   & = & -A_u + C_u \\
{}   & = & -2 A_u +\max(f) \\
{}   & = & \max(f)-2\mu_i(x)
\end{eqnarray*}
By theorem 3(1), each
$D_i(x)$ is a totally $2\pi$-convex, totally geodesic submanifold of
$(M, g)$.

Now we will show that the functions $\mu_i(x)$'s are all constant
functions on $D_{\max}$ in the following
\begin{lemma}
\begin{enumerate}
\item For each
$i\in\{1, 2, \cdots ,p(x)\}$ the function $\mu_i(x)$ is constant on $\cm$.
\item The function $p(x)$ is constant on $D_{\max}$.
\end{enumerate}
\end{lemma}
\noindent{\bf Proof:} For each $x\in \cm$, $\lpx$ is the largest
eigenvalue of $-\nabla^2f(x)$. Since the geodesics
$\gamma_u$ are all integral curves of
$\nabla f$ for $u\in U_{\lpx}$ they will flow towards
$D_{\min}:=\{y\in M:f(y)=\min (f)\}$ and they will reach $D_{\min}$
at time $\pi$. Hence we must have
$\max(f)-2\lpx=\max(f)-2\mu_{p(y)}$ for $x$, $y\in \cm$. This
proves that $\mu_{p(x)}$ is a constant function on $\cm$.

Now we will prove the following
\begin{sublemma}
For each critical value $\alpha$ of the function $f$, the submanifold
$D_{\alpha}$ coincides with $D_i(x)$ for every $x\in\cm$ and for some
$i\in \{ 1, 2, \cdots , p(x)\}$.
\end{sublemma}
\noindent{\bf Proof:} Let $y\in D_{\alpha}$ and $x\in\cm$. Then
$d(x, y)=d(D_{\alpha}, \cm)=\pi$. Since $D_{\alpha}$ and $\cm$ are
totally geodesic submanifolds of $(M, g)$, any geodesic segment
$\gamma_u$
joining $\gamma_u(0)=x$ and $\gamma_u(\pi)=y$ will meet both $D_{\alpha}$
and $\cm$ orthogonally.

Let $u_{\theta}:=\cos \theta w +\sin \theta v$ be a curve in $U_xM$
such that $w\in E_{\lpx}$, $v\in E_{\lpx}^{\perp}$ and
$u_{\theta_{0}}=u$ for some $\theta_0$. Then
$f(\gamma_{u_{\theta}}(t))=A_{u_{\theta}}\cos t + C_{u_{\theta}}$ and
\begin{eqnarray*}
A_{u_{\theta}} & = & -\nabla^2f(u_{\theta}, u_{\theta}) \\
{} & = & -\cos^2\theta\nabla^2f(w, w)-\sin^2\theta\nabla^2f(v, v)
							-2\sin\theta\cos\theta <\nabla^2f(w), v> \\
{} & = & \cos^2\theta A_w +\sin^2\theta A_v
					+2\sin\theta\cos\theta A_w <w,v> \\
{} & = & \cos^2\theta A_w + \sin^2\theta A_v
\end{eqnarray*}
The third step in the above equation follows from the fact that $w$ is an
eigenvector of $-\nabla^2 f$ and the last step follows since $<w, v>=0$.

Since $\theta_0$ is a critical point of the function
$\theta\mapsto f(\gamma_{u_{\theta}}(\pi))$, we have that
${d\over {d\theta}}\!\mid_{\theta=\theta_0}f(\gamma_{u_{\theta}}(\pi))=0$.
Further, since
\begin{eqnarray*}
f(\gamma_{u_{\theta}}(\pi)) & = & -A_{u_{\theta}} + C_{u_{\theta}} \\
{} & = & -2 A_{u_{\theta}} + \max(f) \\
{} & = & \max(f) -2(\cos^2\theta A_w +\sin^2\theta A_v)
\end{eqnarray*}
we have that
\begin{eqnarray*}
0 & = & {d\over {d\theta}}\!\mid_{\theta=\theta_0}
										f(\gamma_{u_{\theta}}(\pi)) \\
{} & = & 2\sin 2\theta_0 (A_v-A_w)
\end{eqnarray*}
Since $-\nabla^2f(w, w)=A_w$ is the largest eigenvalue, we have that
$A_v-A_w\neq 0$. Hence $\sin 2\theta_0=0$. i.e.,$\theta_0={{\pi}\over 2}$.
This shows that $u\in E_{\lpx}^{\perp}$. Now from
$\min -\max$ principle it follows that $\alpha =\mu_i(x)$
for all $x\in\cm$ and for some $i\in\{1, 2, \cdots, p(x)\}$.
This completes the proof of the sublemma.

{}From the sublemma it follows that
\begin{enumerate}
\item Each eigenvalue $\mu_i(x)$ is a constant
function on $\cm$
and the number of distinct eigenvalues of $-\nabla^2f$ are constant on
$\cm$. Hence $p(x)=p$, a constant independent of the point $x\in\cm$.
\item The only critical values of the function $f$ are
$\max(f) -2\mu_i$ where $1\leq i\leq p$ and $\mu_i$'s are the
eigenvalues of $-\nabla^2f$ on $\cm$.
\end{enumerate}
This proves the lemma.

Therefore $D_i :=\{y\in M : f(y)=\max(f)-2\mu_i ~\&~ \nabla f(y)=0\}$
are the only critical submanifolds of the function $f$ with critical
values $\alpha_i=\max(f)-2\mu_i$.
Note that $D_1=\cm$ and $D_p=D_{\min}$.

For each eigenvalue $\mu$ of $-\nabla^2f$ on $D_i$, we denote by
$E_{\mu}$, the $\mu$-eigensubspace of $-\nabla^2f$ and by
$S_{\mu}(0, \pi)$ the sphere of radius $\pi$ in $E_{\mu}$.

Then  we have the following
\begin{propn}
\begin{enumerate}
\item The eigenvalues of $-\nabla^2f$ on $D_i$ are
$\{\mu_{ij}:=\mu_j-\mu_i:1\leq j\leq p\}$.
\item For each $x\in D_i$, the map
			$$
				\exp_x\!\mid_{S_{\mu_{ij}}(0, \pi)} :S_{\mu_{ij}}(0, \pi)\to D_j
			$$
is a fibration for each $j\neq i$.
\item Each $D_j$ is either
\begin{enumerate}
\item an integral cohomology CROSS and the degree of the generator
of $H^*(D_j, Z\!\!Z)$ is $k=2$, $4$, $8$ or $n$, or
\item diffeomorphic
to $I\!\!RI\!\!P^{d_{ij}-1}$ where $d_{ij}=dimE_{\mu_{ij}}$.
\end{enumerate}
\end{enumerate}
\end{propn}
\noindent{\bf Proof:} Let $E_{\mu}(x)$ be an eigensubspace of
$-\nabla^2f(x)$ on $D_i$ for $\mu\neq 0$. Since $\exp_x(S_{\mu}(0, \pi))$
is a critical submanifold of the function $f$ it must be one of the
$D_j$'s for $j\neq i$.

Let $u\in E_{\mu}(x)$ be a unit vector. Then
$f(\gamma_u(0))=A_u+C_u=\max(f)-2\mu_i$ and
$f(\gamma_u(\pi))=-A_u+C_u=\max(f)-2\mu_j$. Since
$\mu=-\nabla^2f(u, u)=A_u$, we have that
\begin{eqnarray*}
\max(f)-2\mu_j & = & -A_u + C_u \\
{}             & = & -2 A_u + A_u + C_u \\
{}             & = & -2\mu + \max(f)-2\mu_i
\end{eqnarray*}
Therefore $\mu=\mu_j-\mu_i$. This proves the first part of the
lemma.

Let $x\in C_i$ and
$D_{\mu_{ij}}(0, \pi):=\{v\in
E_{\mu_{ij}}(x):\parallel\!v\!\parallel\leq\pi\}$.
Let $M_{ij}(x):=\exp_x(D_{\mu_{ij}}(0, \pi))$.
Then each $M_{ij}(x)$ is a submanifold of $M$ and
$M_{ij}(x)$ is also a Blaschke manifold at $x$ with totally geodesic
cut-locus $D_j$. Hence it follows from
\cite{O} and \cite{NS} that
$\exp_x\!\mid_{S_{\mu_{ij}}(0, \pi)}:S_{\mu_{ij}}(0, \pi)\to D_j$
is a fibration for $j\neq i$.

Since $(M, g)$ is a $P$-manifold, the index of geodesics is a constant,
say $(k-1)$. From the fact that each $M_{ij}(x)$ is a Blaschke manifold
at $x$, it follows that $k-1\in\{0, 1, 3, 7, n-1\}$ and the dimension of
the fibres is $k-1$ for all these fibrations. We note that this
fact can also be verified combinatorially.

If $k-1$ is positive then each $M_{ij}(x)$ is simply connected
integral cohomology CROSS and the degree of the generator of
$H^*(M_{ij}(x), Z\!\!Z)$ is $k$. Since $D_j$ is the cut-locus of
$x$ in $M_{ij}(x)$, it follows that $D_j$ is also a simply connected
inetgral cohomology CROSS and the degree of the generator of
$H^*(D_j, Z\!\!Z)$ is $k$.

Since each $M_{ij}(x)$ is a Blaschke manifold at $x\in D_i$, for any
point  $y\in D_j=Cut(x)$ in $M_{ij}(x)$, the order of conjugacy is
atleast
$dim\Lambda(x, y)$ along any geodesic $\gamma$ with recpect to $x$ ;
here $\Lambda(x, y)$ is the link between $x$ and $y$.(See
\cite{BE}). Therefore
if the index $k-1$ is zero, then $dim\Lambda(x, y)=0$ for every
$y\in D_j$. Now since $x$ is not conjugate along any geodesic $\gamma$
from $y$ to $x$ and $\Lambda(x, y)$ has only two elements, we see that
$\exp_x\!\mid_{S_{\mu_{ij}}(0, \pi)}:S_{\mu_{ij}}(0, \pi)\to D_j$ is a two
sheeted covering. This proves that $D_j$ is diffeomrphic to
$I\!\!RI\!\!P^{d_{ij}-1}$.
This proves 3, 4 and the proposition completely.
\section{Proof of theorem 2}
Let $\lambda$ be an eigenvalue of $\Delta$ with an eigenfunction $f$ such
that $f(\gamma_u(t)) = A_u\cos t + B_u\sin t + C_u$ for $u\in UM$.
We know from theorem 3 that the function has only finitely many critical
values say $\{\alpha_i:1\leq i\leq p\}$.
Let $\cm = D_1, D_2, \cdots , D_p=D_{\min}$ be the critical
submanifolds of the function $f$ with critical values $\alpha_i$.

Let $\mu_p>\mu_{p-1}>\cdots >\mu_2>\mu_1=0$ be the eigenvalues of
$-\nabla^2f$ on $\cm$. We saw in proposition 1 that for each $x\in\cm$,
the map $\exp_x\!\mid_{S_{\mu_{ij}}(0, \pi)}: S_{\mu_{j}}(0, \pi)\to D_j$
is fibration with fibres of dimension $k-1$. Therefore we can write $dim
E_{\mu_{j}}=kr_j$ for some non-negative integer $r_j\in\{1, 2, \cdots ,
n\}$. Hence $dimD_j=k(r_j-1)$.

We also saw in proposition 1 that the eigenvalues of $-\nabla^2f$ on
$D_i$ are $\{\mu_{ij}:\mu_j-\mu_i:1\leq j\leq p\}$ and
$\exp\!\mid_{S_{\mu_{ij}}(0, \pi)}:S_{\mu_{ij}}(0,\pi)\to D_j$ is a
fibration for $j\neq i$. In particular
$\exp\!\mid_{S_{-\mu_{i}}(0, \pi)}: S_{-\mu_{i}}(0, \pi)\to D_{\max}$ is a
fibration. Hence $dim E_{\mu_{ij}}=dim E_{\mu_{j}}=kr_j$ and
$dim E_{-\mu_{i}}=dim\cm +k=k(r_1+1)$.

Now we will compute the Laplacian of the function $f$ on $D_i$'s.

Since $f$ is an eigenfunction of $\Delta$ with eigenvalue $\lambda$,
for each $x\in\cm$
\begin{eqnarray*}
\lambda\max(f) & = & \Delta f(x) \\
{} & = & Tr(-\nabla^2f(x)) \\
{} & = & k\sum_{i=1}^pr_i\mu_i
\end{eqnarray*}
and for each $y\in D_j$
\begin{eqnarray*}
\lambda\alpha_j & = & \Delta f(y)
\end{eqnarray*}
But we know that $\alpha_j=\max(f)-2\mu_j$. Therefore
\begin{eqnarray*}
\lambda(\max(f)-2\mu_j) & = &k(r_1+1)(\mu_1-\mu_j)
														+k\sum_{i\geq 2}r_i(\mu_i-\mu_j) \\
{} & = & -k\mu_j+k\sum_{i=1}^pr_i(\mu_i-\mu_j) \\
{} & = & -k\mu_j+k\sum_ir_i\mu_i-k\mu_j\sum_ir_i \\
{} & = & -{{k(1+\sum_ir_i)}\over 2}\mu_j+\lambda\max(f)
\end{eqnarray*}
This proves that
								$$
									\lambda = {{k(m+1)}\over 2}
								$$
where $m=\sum_ir_i$.

We know from Bott-Samelson theorem for $P$-manifolds that
$H^*(M, I\!\!\!Q)$ has exactly one generator.(See \cite{BB}, \cite{BE}).
{}From proposition 1 it follows that the degree of the generator is $k$.
Therefore
$\lambda={{k(m+1)}\over 2}=\lambda_1({\olimm})$ where $\olimm$ is a CROSS
of dimension $km$ with sectional curvarture
${1\over 4}\leq K_{\olimm}\leq 1$ and
$H^*(M, I\!\!\!Q)=H^*(\olimm , I\!\!\!Q)$. This proves theorem 2.
\section{Proof of theorem 1}
By hypothesis $Ric_M\geq l$ and
$\lambda_1 = {1\over 3}(2l+n+2)$. Hence for any first
eigenfunction $f$ we have that
$f(\gamma_u(t))=A_u\cos t+B_u\sin t+C_u$ for $u\in UM$.(See \cite{AR}).

\noindent{\bf Proof of 1a:}
It follows from theorem 2 that $\lambda_1= {{k(m+1)}\over 2}$.
Since $\lambda_1$ is also equal to $= {1\over 3}(2l+n+2)$, we
get that $l= {{k(m-1)}\over 4} + (k-1)=Ric_{\olimm}$.
Again from the proof of theorem 2 it follows that
$H^*(M, I\!\!\!Q)=H^*(\olimm, I\!\!\!Q)$ and also that
$H^*(M, Z\!\!Z_2)=H^*(\olimm, Z\!\!Z_2)$.
This completes the proof of 1a.

\noindent{\bf Proof of 1b:} If $k\geq 2$ then each $D_i$ is a simply
connected submanifold of $M$ of dimension less than or equal to $km-k$.
If one of the $D_i$ is of dimension $km-k$, then there are only two
critical submanifolds $\cm$ and $D_{\min}$ of the function $f$ and one of
them is a point. Let us assume that $D_{\min}=\{p\}$. Then $M$ is a
Blaschke manifold at $p$ with simply connected cut-locus $\cm$. Hence $M$
is simply connected.

Now we assume that all the critical submanifolds are of dimension less
than or equal to $km-2k$. For each $i$, the critical submanifold
$D_i$ is contained in $M\setminus\cup_{j\neq i}D_j$ and $D_i$ is a strong
deformation retract of $M\setminus\cup_{j\neq i}D_j$. Since the
codimension of each $D_j$ is greater than or equal to $4$ we have that
$\pi_1(M)=\pi_1(M\setminus\cup_{j\neq i}D_j)$. Therefore
$\pi_1(D_i)=\pi_1(M\setminus\cup_{j\neq i}D_j)=\pi_1(M)$. This proves that
$M$ is simply connected.

We have seen in proposition 1 that each $D_j$ is a simply conected
integral cohomlogy CROSS and the degree of the generator of
$H^*(D_j, Z\!\!Z)$ is $k$. Further we are attaching only $rk$ dimensional
cells to each of these $D_j$'s. This proves that the integral
cohomology ring of $M$ is same as that of $\olimm$.

\noindent{\bf Remark:} If the integral cohomology ring of $M$
is same as that of the cohomology projective plane then the function can
have atmost three critical submanifolds $\cm$ and $D_{\min}$ and one
saddle. If there are three critical submanifolds then all of them are
points; if there are  only two critical submanifolds $D_{\max}$ and
$D_{\min}$ again one of them is a point.

\noindent{\bf Proof of 2:}
Since $k=1$ we have that $Ric_M\geq {{n-1}\over 4}$.
Let $(\wm, \wg)$
be the universal cover of $(M, g)$ and $\Pi: \wm \to M$ the covering
map. Then, since $Ric_M=Ric_{\wm}$, we have that
$Ric_{\wm}\geq {{n-1}\over 4}$. Hence by Bonnet-Myers theorem
$diam(\wm, \wg)\leq 2\pi$. Now we will show that
$diam(\wm, \wg)\geq 2\pi$. Then from the rigidity of Bonnet-Myers
theorem it will follow that $(\wm, \wg)$ is isometric to $S^n$ with
constant sectional curvature ${1\over 4}$ and our proof will also show that
$M$ is isometric to $\rpn$ with sectional curvature
${1\over 4}$.

Let us fix a point $x_0\in\cm$ and let $D(0, \pi)\subseteq T_{x_{0}}M$ be
the disk of radius $\pi$ in $T_{x_{0}}M$. Then $\exp_{x_{0}}:D(0,
\pi)\to M$ is a smooth on-to map. We identify the antipodal
points in the boundary $S(0, \pi)$ of $D(0, \pi)$ and denote it
by $\rpn:=D(0, \pi)/<(u, -u):u\in S(0, \pi)>$, the quotient
space of this identification.

Since $M$ does
not have conjugate points along the geodesics $\gamma$ up to length
$\pi$, $\exp_{x_{0}}: \rpn\to M$ is a smooth on-to map of maximal
rank. Therefore $\exp_{x_{0}}:\rpn\to M$ is a covering and the map
$\Pi: \wm \to M$ factors through $\rpn$. Since
$\exp_{x_{0}}:\rpn\to M$ is a covering, we know that the map
$(\exp_{x_{0}})_*:\pi_1(\rpn)\to \pi_1(M)$ is injective. since the geodesic
loops joining different critical levels are non-trivial in $\rpn$
they are non-trivial in $M$ also. Now let us fix one such geodesic loop
$\gamma$ of length $2\pi$ in $\pi_1(M)$. Let $\wga$ be the lift of
$\gamma$ in $\wm$ with $\wga(0)=x$ and $\wga(2\pi)=y$. We claim that
$d(x, y)=2\pi$. Suppose not. Then there exists a geodesic segement
$\wgaa$ of length less than $2\pi$ joining $x$ and $y$. Since $\wm$ is
simply connected there exists a homotopy between $\wga$ and $\wgaa$
fixing the end points $x$ and $y$. This homotopy will go down to $M$ to
give a homotopy between $\gamma$ and $\gamma_1=\Pi(\wgaa)$. Since the
geodesic $\gamma_1$ is based at a critical point and of length less than
$2\pi$, it can be homotoped to the base point along the integral curves
of $\nabla f$. Hence $\gamma$ is also
trivial in $\pi_1(M)$, a contradiction. Therefore $diam(\wm, \wg)\geq
2\pi$.
This proves that $(\wm, \wg)$ is isometric to $S^n$ with constant
sectional curvature ${1\over 4}$ and $(M, g)$ is isometric to $\rpn$ with
constant sectional curvature ${1\over 4}$.
This completes the proof.

\noindent{\bf Proof of 3:}
Since $k=n$ we have that $Ric_M\geq n-1$. Further we also know $diam (M,
g)=\pi$.
Hence it follows from Bonnet-Myers theorem that $(M, g)$ is
isometric to $S^n$ with constant sectional curvature $1$.

\noindent{\bf Remark:}
We note that, since $k=n$, the function $f$ has only two critical
points, namely the maxima and minima and they are non-degenerate. Hence
$f$ does not have saddle points.
\subsection{Proof of 4}
In this subsection we will assume that $\max(f)$ and $\min(f)$ are the
only critical values of the function $f$. Hence $D_{\max}$ and
$D_{\min}$ are
the only critical submanifolds of the function $f$ in $(M, g)$.
Therefore $-\nabla^2f$ has only two eigenvalues on $D_{\max}$. By
normalising the function $f$, we may assume that these two eigenvalues are
$1$ and $0$. Hence we can write $f(\gamma_u(t))=\cos t + C$ for
$u\in (UD_{\max})^{\perp}$, the unit normal bundle of $D_{\max}$.
Since the integral curves of
$\nabla f$ are geodesics, the normal geodesic spheres around
$D_{\max}$ are
level sets of the function $f$.

Now we get bounds for $\nabla^2f(u, u)$ for every $u\in UM$.

Let $S(t)$
be the geodesic sphere of radius $t$ around $D_{\max}$. Then
$f(x) = \cos t + C$ for $x\in S(t)$ and
$f(\gamma_u(t))=A_u\cos t +B_u\sin t +C_u$
for $u\in U_xM$. Then $\gamma_u(0)\in S(t)$ and
$\gamma_u(\pi)\in S(t_1)$ for some $t_1$ such that $0\leq t_1\leq \pi$.
Since $A_u+C_u=\cos t+C$ and $-A_u+C_u=\cos t_1+C$, we have that
$A_u={1\over 2}(\cos t-\cos t_1)$. Therefore
\begin{eqnarray*}
-\nabla^2f(u, u) & = & A_u \\
{}               & = & {1\over 2} (\cos t -\cos t_1) \\
\end{eqnarray*}
Since $-1\leq \cos t_1\leq 1$, we get that
				$$
					{{1-\cos t}\over 2}\geq \nabla^2f(u, u)
														 \geq {-{{1+\cos t}\over 2}}
				$$
Having got these bounds for $\nabla^2f$, we define two eigensubbundles
of $\nabla^2f$
\begin{eqnarray*}
\eeu & := &
				\{E\in TM : \nabla^2f(E)= {{1-\cos t}\over 2}\} \\
\euu& := &
				\{E\in TM : \nabla^2f(E)= {-{{1+\cos t}\over 2}}\}
\end{eqnarray*}
Then we have the follwing
\begin{lemma}
\begin{enumerate}
\item The eigensubbundles $\eeu$ and $\euu$ of $\nabla^2f$ are parallel
along $\nabla f$.
\item $\eeu$ and $\euu$ are eigensubbundles of
$R(., \nabla f)\nabla f$ with eigenvalue ${1\over 4}\minnu^2$.
\end{enumerate}
\end{lemma}
\noindent{\bf Proof:} Let $x\in\cm$ and $\gamma$ be a geodesic starting
at $x$ such that $\gamma'(0)\in (U\cm)^{\perp}$.
Let $J$ be a Jacobi field along $\gamma$ describing the
variation of the geodesic $\gamma$ such that $J(0)\in TD_{\max}$
and $J(\pi)=0$. We normalise $J$ such that
$\parallel\!J'(\pi)\!\parallel=1$. Then, since $J$ is a Jacobi field,
$[J, \gamma'(t)]=0$ along the geodesic $\gamma$. Further, since
$\gamma'(t) ={{\nabla f}\over \minnu}$, we have that
$\gamma'(t) = -\partial_t$. Hence
\begin{eqnarray*}
-<J', J> & = & {1\over \minnu}<\nabla_J\nabla f, J> \\
{}       & \leq & {{\!\parallel J\!\parallel}^2\over \minnu}
									{{1-\cos t}\over 2} \\
{{<J', J>}\over{\parallel J\parallel}^2} & \geq &
							-{1\over 2}{{\sin{t\over 2}}\over{\cos{t\over 2}}}
\end{eqnarray*}
The function ${{\parallel J\parallel}^2\over{\cos^2{t\over 2}}}$ is
smooth  and non-vanishing on $I\!\!R$. Hence we can take the positive
square root
${{\parallel J\parallel}\over{\mid\cos{t\over 2}\mid}}$ of
${{\parallel J\parallel}^2\over{\cos^2{t\over 2}}}$
which is again smooth.
Since the function $\cos t$ is an even function,
$\cos{t\over 2}$ is positive on $(-\pi, \pi)$.
Therefore from the last step of the above equation it follows that
\begin{eqnarray*}
\dt \log({{\parallel J\parallel}\over{\cos{t\over 2}}}) & \geq & 0
\end{eqnarray*}
on $(-\pi, \pi)$.
Now since $(M, g)$ is a $\ptp$-manifold, we have that
$J(t)= J(t+2\pi)$. Hence
${{\parallel J \parallel}\over{\cos{t\over 2}}}\mid_{t=-\pi}=
{{\parallel J \parallel}\over{\cos{t\over 2}}}\mid_{t=\pi} = 2$.
This proves that
${{\parallel J \parallel}\over{\cos{t\over 2}}}=2$ for
$t\in [-\pi, \pi]$ and equality must hold everywhere in the above
inequalities. This proves that $J$ is an eigenvectorfield of
$\nabla^2f$ with eigenvalue ${{1-\cos t}\over 2}$.
Since
$\parallel J\parallel = 2\cos{t\over 2}$, we can write
$J(t) = 2\cos{t\over 2} E(t)$ where $E(t)\in\eeu$ is a unit vector field
along $\gamma$. Since $J$ is a Jacobi field along $\gamma$
\begin{eqnarray*}
J' & = & \nabla_J\gamma' \\
{} & = & {1\over\minnu}\nabla_J\nabla f \\
{} & = & {{1-\cos t}\over 2}{1\over\minnu}J \\
{} & = & {{1-\cos t}\over 2}{1\over\minnu}2\cos{t\over 2}E
\end{eqnarray*}
on the other hand
$J'=\sin{t\over 2}E+\cos{t\over 2}E'$. This shows that $E'$ is along the
direction of the vector field $E$. Since $E$ is a unit vector field along
$\gamma$, $E'\perp E$. Therefore $E'=0$ along $\gamma$. This proves that
$\eeu$ is parallel along $\nabu$.

Now by a similar argument we can show that the eigensubbundle $\euu$ is
also parallel along $\nabu$ by using the inequality that
$\nabla^2f(u, u)\leq -{{1+\cos t}\over 2}$. (For a proof see \cite{RG}).
This completes the proof of the lemma 3(1).

Now we set out to prove lemma 3(2). Let $E\in \eeu$ be a unit vector at
$t=0$ and $J$ be a Jacobi field describing the variation of a normal
geodesic $\gamma$ starting $\cm$,
such that $J(0)=2E$. Then from what we have seen above
$J(t)=2\cos{t\over 2} E(t)$; $E(t)$ parallel along $\gamma$. Therefore
\begin{eqnarray*}
R(J, \gamma')\gamma' & = & - J'' \\
{}                   & = & {1\over 4}J
\end{eqnarray*}
and this proves that $\eeu$ is eigensubbundle of $R(., \nabu)\nabu$ with
eigenvalue ${1\over 4}\minnu^2$ along $\nabu$. The same
arguments will prove that $\euu$ is also an eigensubbundle of
$R(., \nabu)\nabu$ with eigenvalue ${1\over 4}\minnu^2$.
This completes the proof of the lemma.

Let $dim\cm=ka$ and $dim D_{\min}=kb$ for some non-negative integers $a$
and $b$. Then $dim\eeu=ka$ and $dim\euu=kb=k(m-a+1)$.

Let $E_{-\cos t} := (\eeu \oplus\euu)^{\perp}$ be the orthogonal
complement of $\eeu\oplus\euu$ in $TM$. Then we have the following
\begin{lemma}
$\eu$ is an eigensubbundle of
\begin{enumerate}
\item $\nabla^2f$ with eigenvalue $-\cos t$
\item $R(., \nabu)\nabu$ with eigenvalue $\minnu^2$
\end{enumerate}
\end{lemma}
\noindent{\bf Proof:}
First we note that $dim(\eeu \oplus\euu)=k(m-1)$. Therefore the
dimension of $\eu$ is $k$. Let us choose an orthonormal basis
$E_2, E_3, \cdots , E_k$ of $\eu$,
$E_{k+1}, E_{k+2}, \cdots , E_{ka}$ of $\eeu$ and
$E_{ka+1}, E_{ka+2}, \cdots , E_{km}$ of $\euu$. Then
\begin{eqnarray*}
\sum_{i=2}^k <R(E_i, \nabu)\nabu, E_i > & = & Ric_M(\nabu , \nabu)
												-\sum_{j=k+1}^{kn}<R(E_j, \nabu)\nabu, E_j> \\
{} & = & \big[{{k(m-1)}\over 4} + (k-1)\big]\minnu^2
					-{{k(m-1)}\over 4}\minnu^2 \\
{} & = & (k-1)\minnu^2
\end{eqnarray*}
Now, for $2\leq i\leq k$, we define the vector fields $W_i=\sin t E_i(t)$,
where each $E_i$ is a parallel vector field along $\gamma$ such that
$E_i(0)=E_i$. Then from the Index lemma, it follows that
$0\leq I(W_i, W_i) =\int_0^{\pi}(<W_i', W_i'>-R(W_i, \gamma')\gamma',W_i>)$.
Therefore
\begin{eqnarray*}
0 & \leq & \sum_{i=2}^kI(W_i, W_i) \\
{} & = & \sum_{i=2}^k\int_0^{\pi}\cos^2t<E_i, E_i>-
						\sin^2t K(E_i,\gamma') \\
{} & = & (k-1)\int_0^{\pi}(\cos^2t-\sin^2t) \\
{} & = & 0
\end{eqnarray*}
Hence $W_i=\sin t E_i(t)$ are Jacobi fields along $\gamma$ for
$2\leq i\leq k$. Now it can be easily verified that $\eu$ is an
eigensubbundle of $\nabla^2f$ with eigenvalue $-\cos t$ and also
an eigensubbundle of
$R(., \gamma')\gamma'$ with eigenvalue $1$. This completes the proof of
the lemma.

\noindent{\bf An interesting Remark:}
When $k=2$, we don't need the condition on $Ric_M$ to
show that $\eu$ is an eigensubbundle of $\nabla^2f$ with
eigenvalue $-\cos t$ and also an eigensubbundle of
$R(.,\gamma')\gamma'$ with eigenvalue $1$. We give the proof below.

Let $x\in\cm$. Then
\begin{eqnarray*}
\Delta f(x) & = & {{k(m+1)}\over 2}f(x) \\
{}          & = & {{k(m+1)}\over 2}(1+C)
\end{eqnarray*}
Therefore
\begin{eqnarray*}
{{k(m+1)}\over 2}(1+C) & = & Tr(-\nabla^2f(x)) \\
{}                     & = & -Tr(\nabla^2f(x)\mid_{\euu})
														 -Tr(\nabla^2f(x)\mid_{\eu})\\
{}                     & = & k(m-a)
\end{eqnarray*}
Hence $C={{m-(2a+1)}\over{m+1}}$.

Now let $p\in M$. Then $f(p)=\cos t+C$ for some $t$ and
\begin{eqnarray*}
{{k(m+1)}\over 2}\big[\cos t +C] & = & Tr(-\nabla^2f(p))\\
{} & = & -\mu_1-\mu_2-Tr(\nabla^2f(p)\mid_{\euu}) \\
{} & {} &  -Tr(\nabla^2f(p)\mid_{\eeu}) \\
{} & = & \cos t -\mu_2- ka({{1-\cos t}\over 2}) \\
{} & {} & 	+k(m-(a+1))({{1+\cos t}\over 2})
\end{eqnarray*}
Hence by substituting the value ${{m-(2a+1)}\over{m+1}}$ for $C$ we
get that $\mu_2=-\cos t$. This completes the proof.

An important consequence of lemma 3 is that, for each $x\in\cm$ the map
$\exp_x: S(0, \pi)\to D_{\min}$ and for each $y\in D_{\min}$ the map
$\exp_y: S(0, \pi)\to\cm$ are great sphere fibrations; here $S(0, \pi)$
denotes the normal sphere of radius $\pi$ at the corresponding points.
Now we prove the following
\begin{lemma}
For every $x\in \cm$, the map
				$$
						\exp_x : S(0, \pi)\to D_{\min}
				$$
and for every $x\in D_{\min}$, the map
				$$
						\exp_x : S(0, \pi)\to \cm
				$$
are congruent to Hopf fibrations.
\end{lemma}
\noindent{\bf Proof:} See \cite{GG} and \cite{RG}.

\noindent{\bf Proof of 4:}
Let us fix a $I\!\!P^a(k)\subseteq I\!\!P^m(k)$. We denote by
$(T\cm)^{\perp}$, the normal bundle of $\cm$ and by
$(TI\!\!P^a(k))^{\perp}$, the normal bundle of $I\!\!P^a(k)$. Since
the map $\exp_x : S(0, \pi)\to D_{\min}$ is congruent to Hopf
fibration for each $x\in\cm$ there is a fibre preserving isometry
$I : (T\cm)^{\perp}\to (TI\!\!P^a(k))^{\perp}$. Using this isometry we
define a map
				$$
					\Phi : M\setminus D_{\min}\to I\!\!P^n(k)
				$$
as follows: For every $q\in M\setminus D_{\min}$ there is a unique
$x\in\cm$ and a unique geodesic segement joining $x$ and $q$ and
we define \mbox{$\Phi(q):=\exp\circ I\circ\exp_m^{-1}(q)$.}
This map carries the
geodesics orthogonal to $\cm$ to geodesics orthogonal to $I\!\!P^a(k)$
and matches the normal geodesic spheres around $\cm$.
To complete the proof we only have
to show that $d\Phi$ preserves the length of the Jacobi fields along
these normal geodesics. This follows from \cite{RG}.
Hence the proof of theorem1(4).

\noindent{\em Akhil Ranjan, Department of Mathematics, \\
Indian Institute of Technology,
Bombay- 400 076,  India. \\
e-mail: aranjan@ganit.math.iitb.ernet.in
\vskip.5cm
\noindent G.~Santhanam, School of Mathematics, \\
Tata Institute of Fundamental Research,
Bombay- 400 005,  India. \\
e-mail: santhana@tifrvax.tifr.res.in, santhana@math.tifr.res.in
\vskip.75cm
}
\end{document}